# Computer modeling of exponential and logarithmic functions of generalized quaternions in symbolic computation system

Yakiv O. Kalinovsky[1], Yuliya E. Boyarinova[2,1], Dmitry V. Lande[1,2], Alina S. Sukalo[1]

[1]IPRI NAS of Ukraine, Kyiv, Shpaka str. 2, 03113, Ukraine
[2]NTUU "KPI", Kyiv, Peremogy av. 37, 03056, Ukraine,

*The paper deals with the process of mathematical modeling representations of exponential and logarithmic functions hypercomplex number system of generalized quaternions via determining a linear differential equation with hypercomplex coefficients. Simulation is performed using symbolic computation system Maple. Some properties of these concepts and their relation to the exponential representations of specific non-commutative hypercomplex number systems of dimension four.*

**Introduction**

Hypercomplex number systems(HNS) and, in particular, quaternions are widely used in solving many scientific problems and equipment in the mechanics of solids for describing the rotation in space, navigation problems, orientation and motion control, in computer animation, deformation study elastic structures, filtration colored imaging, cryptography, data constructing fast algorithms and many others. Research in this field are currently in [1-7].

Generalized quaternions were introduced K. Gödel in 1949 to represent the space-time groups. In [8] he presented solutions of Einstein's gravitational field equations by using generalized quaternions.

Investigation of generalized quaternions dedicated work of many other authors, for example, [9-15].

In [16] generalized quaternions proposed for increasing the resistance of the electronic digital signature in the RSA algorithm. This requires the creation of mathematical models of various representations of nonlinearities of the generalized quaternion, and, above all, of the exponential function and its inverse function - logarithm of the generalized quaternion.

However, in these studies did not under consideration of constructing representations of various nonlinear functions of generalized quaternions, and, above all, exponential and logarithmic functions.

Generalized quaternion is:
$$q = a_1 e_1 + a_2 e_2 + a_3 e_3 + a_4 e_4, \qquad (1)$$

where $a_i$ - real number, $e_i$, $i = 1,...,4$ - basic elements that match the following table Cayley:

| $H_{\alpha\beta}$ | $e_1$ | $e_2$ | $e_3$ | $e_4$ |
|---|---|---|---|---|
| $e_1$ | $e_1$ | $e_2$ | $e_3$ | $e_4$ |
| $e_2$ | $e_2$ | $-\alpha e_1$ | $e_4$ | $-\alpha e_3$ |
| $e_3$ | $e_3$ | $-e_4$ | $-\beta e_1$ | $\beta e_2$ |
| $e_4$ | $e_4$ | $\alpha e_3$ | $-\beta e_2$ | $-\alpha\beta e_1$ |

(2)

where $\alpha$, $\beta \in R$.

## 1. Formulation of the problem

Determination of exponential and logarithmic functions is based on the idea of infinite power series, as suggested by the founder of the theory of hypercomplex numbers W.R.Hamilton [17]. This creates great difficulties in the use of these functions in the simulation.

The paper deals with the construction of such representations exponential and logarithmic functions of the generalized quaternion in which these infinite series of are collapsed to final form. This enables them to be successfully used in the computer simulation.

Due to the complexity of the implementation of the necessary reforms of the problem it was only possible using symbolic computation system Maple.

As shown by experimental research [18], the use of representations of exponential and logarithmic functions significantly reduces the number of calculations.

## 2. The method of constructing representations of exponential hyper variable via the associated system of linear differential equations

Later hypercomplex numbers will be denoted by upper case Latin letters:

$$X = \sum_{i=1}^{n} x_i e_i ; \quad M = \sum_{i=1}^{n} m_i e_i ,$$  (3)

a column vectors composed of components hypercomplex numbers - upper case Latin letters with a dash

$$\overline{X} = (x_1,..., x_n)^T , \quad \overline{M} = (m_1,..., m_n)^T .$$  (4)

Then the main provisions of the proposed method of construction of representations of exponential hyper variable via the associated system of linear differential equations, which are presented in [19 - 23], are as follows.

Presentation of the exhibitors in the hypercomplex number system $\Gamma(e,n)$, the number of which $M \in \Gamma(e,n)$ will be denoted $Exp(M)$, there is a particular solution of the linear differential equation hypercomplex

$$\dot{X} = MX$$  (5)

with the initial condition

$$Exp(0) = \varepsilon ,$$  (6)

where $\varepsilon$ - a single element of the system $\Gamma(e,n)$.

Equation (5) is called the defining.
To construct solutions hyper linear differential equation (5) it must be represented in vector-matrix form:

$$\dot{\overline{X}} = \overline{MX} .$$  (7)

In this case a column vector on the left side is of the form:

$$\dot{\overline{X}} = (\dot{x}_1,...,\dot{x}_n)^T .$$  (8)

Since hypercomplex number is the product of two hyper complex numbers $MX$, it looks like:

$$MX = \sum_{i=1}^{n}\sum_{j=1}^{n}\sum_{k=1}^{n} m_i x_j \gamma_{ij}^k e_k, \tag{9}$$

where $\gamma_{ij}^k$ - structure constants of HNS. Therefore, a column vector $\overline{MX}$ is of the form:

$$\overline{MX} = \left( \sum_{i=1}^{n}\sum_{j=1}^{n} m_i x_j \gamma_{ij}^1, ..., \sum_{i=1}^{n}\sum_{j=1}^{n} m_i x_j \gamma_{ij}^n \right)^T. \tag{10}$$

Column vector $\overline{MX}$ can be represented as a matrix product the size $n \times n$ of a matrix M by a vector-column $\overline{X}$:

$$\overline{MX} = M\overline{X}. \tag{11}$$

Wherein

$$M_{jk} = \sum_{i=1}^{n} m_i \gamma_{ij}^k. \tag{12}$$

Then hypercomplex equation (7) will turn into a system of $n$ linear differential equations of first order, which is called the associated system of linear differential equations

$$\dot{\overline{X}} = M\overline{X}. \tag{13}$$

Next, you need to find the characteristic values $\lambda_1, ..., \lambda_n$ of the matrix M, that is, to solve the characteristic equation

$$M - \lambda E = 0. \tag{14}$$

Thus, the characteristic values of the components will depend on the part of hypercomplex number $M$.

Then construct the general solution depends on $n^2$ arbitrary constants of which are $n^2 - n$ linearly dependent on the $n$ available variables. Excluding dependent constant, we can obtain the general solution of the equation (7), depending on $n$ arbitrary constants—$\overline{X}(t, C_1, ..., C_n)$. The values of the arbitrary constants are set by the initial conditions (6). The components of the column vector $\overline{X}$ of the solution and will be exhibitors from hypercomplex numbers $M$

$$Exp(M) = \sum_{i=1}^{n} \overline{x}_i e_i. \tag{15}$$

### 3. Construction of exhibitors presenting a generalized quaternion

In accordance with (9) vector-matrix system of equations (7) takes the form:

$$\begin{aligned}
\dot{x}_1 &= m_1 x_1 - \alpha m_2 x_2 - \beta m_3 x_3 - \alpha\beta m_4 x_4 \\
\dot{x}_2 &= m_2 x_1 + m_1 x_2 - \beta m_4 x_3 + \beta m_3 x_4 \\
\dot{x}_3 &= m_3 x_1 + \alpha m_4 x_2 + m_1 x_3 - \alpha m_2 x_4 \\
\dot{x}_4 &= m_4 x_1 - m_3 x_2 + m_2 x_3 + m_1 x_4
\end{aligned} \tag{16}$$

a equation (14) becomes the equation for the parameter $\lambda$.

$$\begin{vmatrix} m_1 - \lambda & -\alpha m_2 & -\beta m_3 & -\alpha\beta m_4 \\ m_2 & m_1 - \lambda & -\beta m_4 & \beta m_3 \\ m_3 & \alpha m_4 & m_1 - \lambda & -\alpha m_2 \\ m_4 & -m_3 & m_2 & m_1 - \lambda \end{vmatrix} = 0. \tag{17}$$

If we introduce the notation:

$$\begin{vmatrix} m_1 - \lambda & -\alpha m_2 \\ m_2 & m_1 - \lambda \end{vmatrix} = D_1, \qquad \begin{vmatrix} m_3 & \alpha m_4 \\ m_4 & -m_3 \end{vmatrix} = D_2, \tag{18}$$

then from (17) in the Laplace Theorem we have:
$$D_1^2 + \beta D_2^2 = 0. \tag{19}$$

The solutions of the equation (19), taking into account the determinants of values (18) are two pairs of double roots:
$$\lambda_{1,2} = m_1 + \sqrt{\overline{m}}, \quad \lambda_{3,4} = -m_1 + \sqrt{\overline{m}}, \tag{20}$$

where
$$\overline{m} = -(\alpha m_2^2 + \beta m_3^2 + \alpha\beta m_4^2). \tag{21}$$

Since the sign $\overline{m} = -(\alpha m_2^2 + \beta m_3^2 + \alpha\beta m_4^2)$ not known, the solution of equation (13) must be solved in the form of:
$$x_i = (C_{i1} + tC_{i2})e^{(m_1+\sqrt{\overline{m}})t} + (C_{i3} + tC_{i4})e^{(-m_1+\sqrt{\overline{m}})t}, \quad i = 1,...,4. \tag{22}$$

Of the 16 arbitrary constants $C_{ik}$ are independent only 4. They can be found from the initial conditions of the form (6):

$$x_1(0) = 1, x_2(0) = x_3(0) = x_4(0) = 0. \tag{23}$$

With this in mind, the solution of equation (13) take the form:
$$x_1(t) = \frac{1}{2}(e^{(m_1+\sqrt{\overline{m}})t} + e^{(m_1-\sqrt{\overline{m}})t}), \tag{24}$$

$$x_i(t) = \frac{m_i}{2\sqrt{\overline{m}}}(e^{(m_1+\sqrt{\overline{m}})t} - e^{(m_1-\sqrt{\overline{m}})t}), \quad i = 2,3,4. \tag{25}$$

If $m_i t$ we denote $m_i$, in accordance with (15) we obtain the final expression for the exponential representation of the generalized quaternion:
$$Exp(M) = \frac{1}{2}e^{m_1}[(e^{\sqrt{\overline{m}}} + e^{-\sqrt{\overline{m}}})e_1 + \frac{1}{\sqrt{\overline{m}}}(e^{\sqrt{\overline{m}}} - e^{-\sqrt{\overline{m}}})(m_2 e_2 + m_3 e_3 + m_4 e_4)]. \tag{26}$$

It can be shown that this formula by substituting the respective values of the parameters $\alpha$ and $\beta$ will give an exponent in hypercomplex number system, Cayley table is obtained from the table (2) by substituting in it the same values $\alpha$ and $\beta$. For example, as shown in [[13-15]], and the table when the Cayley's table becomes to table Cayley of ordinary quaternions. Then from (21) we have:
$$\overline{m} = -(m_2^2 + m_3^2 + m_4^2),$$
$$\sqrt{\overline{m}} = i\sqrt{m_2^2 + m_3^2 + m_4^2} = i\sqrt{|\overline{m}|}$$

and using Euler's formula we obtain a representation of the quaternion exponential [20, 22, 24]:
$$Exp(M) = e^{m_1}(\cos\sqrt{|\overline{m}|}e_1 + \frac{\sin\sqrt{|\overline{m}|}}{\sqrt{|\overline{m}|}}(m_2 e_2 + m_3 e_3 + m_4 e_4)). \tag{27}$$

It should be noted that the representation (26) is not satisfied the main characteristic property of exhibitors: exhibitor's arguments amount equal to the product of the exponents of separate arguments:
$$Exp(M_1 + M_2) \neq Exp(M_1) \cdot Exp(M_2).$$

This fact, which can be verified directly by using (26), is due to non-commutative multiplication of generalized quaternions.

However, a remarkable is that kind of representation (26) is invariant with respect to the switching of the factors in determining the equation (5). This is because, as shown in (10), the fact that the solution of the equation (5) can be represented as the sum of a power series

$$X(t) = \sum_{s=0}^{\infty} \frac{(Mt)^s}{s!},$$

which, since associativity of multiplication of generalized quaternions, it has a unique value.

## 4. Building a presentation of the logarithmic function of the generalized quaternion

In accordance with the definition of the logarithmic function as the inverse exponential [25-29] from the expression (26)

$$Ln\{\frac{1}{2}e^{m_1}[(e^{\sqrt{\overline{m}}} + e^{-\sqrt{\overline{m}}})e_1 + \frac{1}{\sqrt{\overline{m}}}(e^{\sqrt{\overline{m}}} - e^{-\sqrt{\overline{m}}})(m_2 e_2 + m_3 e_3 + m_4 e_4)]\} = M. \quad (28)$$

But to represent a logarithmic function of the generalized quaternion needed to the left under the logarithm sign was hypercomplex numbers, and the right - hypercomplex function:

$$Ln \sum_{k=1}^{4} x_k e_k = \sum_{k=1}^{4} f_k(x_1, x_2, x_3, x_4) e_k. \quad (29)$$

To determine the type of function $f_k(x_1, x_2, x_3, x_4)$ to be resolved with respect $m_1, m_2, m_3, m_4$ to a system of four equations:

$$\begin{cases} \dfrac{e^{m_1 + \sqrt{\overline{m}}} + e^{m_1 - \sqrt{\overline{m}}}}{2} = x_1 \\ \dfrac{m_k}{\sqrt{\overline{m}}} \cdot \dfrac{e^{m_1 + \sqrt{\overline{m}}} - e^{m_1 - \sqrt{\overline{m}}}}{2} = x_k, \; k = 2,3,4 \end{cases}. \quad (30)$$

The system (30) is a very difficult structure, since it contains the irrationality of the unknown variables, the signs radicands unknown, and exponential functions, in terms of which are the irrationality. Therefore, its solution is generally difficult. However, this type is greatly simplified by fixing the sign under the radical irrationality. This is because the Euler formula for exponential functions can be reduced to either a trigonometric function or a hyperbolic. It is therefore advisable to consider two cases: $\overline{m} < 0$ and $\overline{m} > 0$.

1) $\overline{m} = -(\alpha m_2^2 + \beta m_3^2 + \alpha\beta m_4^2) < 0.$ \quad (31)

In this case, the use of Euler's formula, and $\sqrt{\overline{m}} = i\sqrt{-\overline{m}}$ the resulting conversion system (30) as follows:

$$\begin{cases} e^{m_1} \cos\sqrt{-\overline{m}} = x_1 \\ e^{m_1} \dfrac{m_k}{\sqrt{-\overline{m}}} \cdot \sin\sqrt{-\overline{m}} = x_k, \; k = 2,3,4 \end{cases}. \quad (32)$$

Multiplying the last three equations, respectively $\alpha$, $\beta$ and $\alpha\beta$, erecting all the equations in the square and adding them, we obtain:

$$e^{2m_1} = x_1^2 + \alpha x_2^2 + \beta x_3^2 + \alpha\beta x_4^2 > 0, \quad (33)$$

from:

$$m_1 = \frac{1}{2}\ln(x_1^2 + \alpha x_2^2 + \beta x_3^2 + \alpha\beta x_4^2). \quad (34)$$

We introduce the notation

$$x_1^2 + \alpha x_2^2 + \beta x_3^2 + \alpha\beta x_4^2 = \overline{x}. \quad (35)$$

Then, from the first equation (32) we obtain:

$$\cos\sqrt{-\overline{m}} = \frac{x_1}{\sqrt{\overline{x}}}, \quad (36)$$

$$\sqrt{-\overline{m}} = Arc\cos\frac{x_1}{\sqrt{\overline{x}}}, \tag{37}$$

$$\sin\sqrt{-\overline{m}} = \sqrt{1-\cos^2\sqrt{\overline{m}}} = \sqrt{1-(\frac{x_1}{\sqrt{\overline{x}}})^2} = \frac{\sqrt{\overline{x}-x_1^2}}{\sqrt{\overline{x}}}. \tag{38}$$

Substitution of (33), (37) and (38) in the last three equations of system (32) produce:

$$m_k = \frac{x_i}{\sqrt{\overline{x}-x_1^2}} Arc\cos\frac{x_1}{\sqrt{\overline{x}}}, \quad k=2,3,4. \tag{39}$$

Transition to (39) on $Arc\cos$ to the main values $arctg$ and substituting in (28) gives the final formula for the logarithmic function in this case:

$$Ln(\sum_{k=1}^{4} x_k e_k) = \ln\sqrt{\overline{x}}\cdot e_1 + \frac{1}{\sqrt{\overline{x}-x_1^2}} arctg\frac{\sqrt{\overline{x}-x_1^2}}{x_1}(\sum_{k=2}^{4} x_k e_k + \pi\sum_{k=2}^{4} n_k e_k), \quad n_k \in Z \tag{40}$$

2) $\overline{m} = -(\alpha m_2^2 + \beta m_3^2 + \alpha\beta m_4^2) > 0$. \quad(41)

In this case, the use of Euler's formula gives the system (30) as follows:

$$\begin{cases} e^{m_1} ch\sqrt{-\overline{m}} = x_1 \\ e^{m_1} \frac{m_k}{\sqrt{-\overline{m}}} \cdot sh\sqrt{-\overline{m}} = x_k, \quad k=2,3,4 \end{cases}. \tag{42}$$

Multiplying the last three equations, respectively $\alpha$, $\beta$ and $\alpha\beta$, erecting all the equations in the square and subtracting the last three of the first, we get:

$$e^{2m_1} = x_1^2 - \alpha x_2^2 - \beta x_3^2 - \alpha\beta x_4^2 > 0, \tag{43}$$

from whence:

$$m_1 = \frac{1}{2}\ln(x_1^2 - \alpha x_2^2 - \beta x_3^2 - \alpha\beta x_4^2). \tag{44}$$

We introduce the notation

$$x_1^2 - \alpha x_2^2 - \beta x_3^2 - \alpha\beta x_4^2 = \overline{x}. \tag{45}$$

Then, from the first equation of system (32) we get:

$$ch\sqrt{-\overline{m}} = \frac{x_1}{\sqrt{\overline{x}}}, \tag{46}$$

$$\sqrt{-\overline{m}} = arch\frac{x_1}{\sqrt{\overline{x}}}, \tag{47}$$

$$\sin\sqrt{-\overline{m}} = \sqrt{ch^2\sqrt{\overline{m}}-1} = \sqrt{(\frac{x_1}{\sqrt{\overline{x}}})^2-1} = \frac{\sqrt{x_1^2-\overline{x}}}{\sqrt{\overline{x}}}. \tag{48}$$

Substitution of (43), (47) and (48) in the last three equations of system (31) produce:

$$m_k = \frac{x_i}{\sqrt{x_1^2-\overline{x}}} arch\frac{x_1}{\sqrt{\overline{x}}}, \quad k=2,3,4. \tag{49}$$

Transition in (48) from $arch$ to $arth$ and a substitution in (27) gives the final formula for the logarithmic function in the second case:

$$Ln(\sum_{k=1}^{4} x_k e_k) = \ln\sqrt{\overline{x}}\cdot e_1 + \frac{1}{\sqrt{x_1^2-\overline{x}}} arth\frac{\sqrt{x_1^2-\overline{x}}}{x_1}\sum_{k=2}^{4} x_k e_k. \tag{50}$$

Using views exponential and logarithmic functions as (27) and (50) resulted in significant reduction number of calculations

**Summary**

The results obtained in this paper presenting exponential and logarithmic functions in the system of generalized quaternions gives the opportunity to get representation for many classes of non-commutative hypercomplex numerical fourth dimension systems. In contrast to existing work to create these views using the method developed by the authors associated with the HNS system of linear differential equations in combination with the use of symbolic computation system. Using presentation exponential and logarithmic functions in the system of generalized quaternions simplifies the analytical construct, and significantly reduces the number of calculations by computer simulation.

**Bibliography**

1. Godel C. An Example of a New Type of Cosmological Solutions of Einstein`s Field Equations of Gravitation // Rev. Mod. Phys. — 1949. — Vol. 21, No. 3. — P. 447—450.
2. Noether E. Hypercomplex Grossen und Darstellungstheorie / Noether E. // Mathematische Zeitschrift. — 1929. — B. 30 P. 641-692.
3. Klipkov S.I. Analysis of generalized matrix representations of associative hypercomplex number systems used in the energy problems / I.S.Klipkov // Data Recording,Storage & Processing. — 2014. — T. 16, № 2. — P. 28–41.
4. Alagos Ya. Split Quaternion Matrices / Ya. Alagos, K. At H. Oral, S. Yuce // Miscolc Mathematical Notes. – 2012. – Vol. 13, No. 2 – P.223–232.
5. Janovska D. Linear equations and the Kronecker product in coquaternions / D. Janovska, G. Opfer // Mitt. Math. Ges. Hamburg. – 2013. Vol. 33. – P. 181–196.
6. Kalinovsky Y. O. Computing Characteristics of One Class of Noncommutative Hypercomplex Number Systems of 4-dimension/ Y. O. Kalinovsky, D. V. Lande, Y. E. Boyarinova, A. S. Turenko. – On-line: http://arxiv.org/ftp/arxiv/papers/1409/1409.3193.pdf.
7. Kalinovsky Y. O. Research links between generalized quaternion and procedures doubling Grassmann-Clifford / Y. O. Kalinovsky, Y. E. Boyarinova, A. S. Turenko // Data Recording,Storage & Processing. — 2015. — T. 17, № 1. — P. 36-45.
8. Mamagami A. B. Some Notes on Matrix of Generalized Quaternion / A. B. Mamagami, M. Jafari // International Research Journal of Applied and Basic Sciences. – 2013. – Vol 7, No. 14. – P. 1164-1171.
9. Kalinovsky Y.OProperties of generalized quaternions and their relation with the procedure of doubling the Grassmann-Clifford / Y. O. Kalinovsky, Y. E. Boyarinova, A. S. Turenko // Electronic Modeling. – 2015, - Vol 37, No. 2.- P.17-26.
10. Sinkov M.V. The finite hypercomplex number systems. Basic theory. Applications / M.V. Sinkov, Y.O.Kalinovsky, Yu.E. Boyarinova — K.: NAS of Ukraine, IPRI, 2010. — 389 p.
11. Kalinovsky Y.O. Methods of computer modeling and calculations using hypercomplex number systems: Dis. PhD : 01.05.02 / Kalinovsky Yakiv Olexandrivich ; IPRI NAS of Ukrane. — K., 2007. — 417 p. —  P. 323–348.
12. Kalinovsky Y.O. Research of isomorphism properties of quadriplex and bicomplex number systems / Y.O.Kalinovsky/Data Recording,Storage & Processing. — 2003. — T. 5, № 1. — P. 69-73.
13. Kalinovsky Y.O. Methods of construction of nonlinear functions in the extended complex numbers / Y.O.Kalinovsky, N.V. Roenko, M.V. Sinkov // Cybernetics and Systems Analysis. — 1996. — № 4. — P. 178–181.

14. Sinkov M.V. Some linear and non-linear operation of generalized complex numbers /Y.O.Kalinovsky, M.V. Sinkov, T.V. Sinkova // Data Recording,Storage & Processing. ─ 2002. ─ Т. 4, № 3. ─ P. 55–61.
15. Brackx F. The Exponential Function of a Quaternion Variable /F. Brackx //Applicable Analysis.-1979.- Vol.19. - P. 265-276.
16. Boyarinova Yu.E. Construction of digital signature algorithm using functions of generalized quaternions / Yu.E. Boyarinova, Y.O.Kalinovsky, A.S.Sukalo// Data Recording,Storage & Processing. ─ 2015. ─ Т. 17, № 3. ─ P. 48–55.
17. Ell, T. A., Bihan, N. L. and Sangwine, S. J. (2014) Quaternion Algebra, in Quaternion Fourier Transforms for Signal and Image Processing, John Wiley & Sons, Inc., Hoboken, NJ, USA.
18. Erik B. Dam Martin Koch Martin Lillholm Quaternions, Interpolation and Animation Technical Report DIKU-TR-98/5 Department of Computer Science University of Copenhagen Denmark July 17, 1998
19. S. Georgiev, J. Morais and W. Sproß New Aspects on Elementary Functions in the Context of Quaternionic Analysis CUBO A Mathematical Journal Vol.14, N¯o01, (93–110). March 2012
20. Kalinovsky Y.O. The logarithmic function of the quaternion / Y.O.Kalinovsky, M.V. Sinkov, T.G. Postnikova, T.V. Sinkova // Data Recording,Storage & Processing. ─ 2002. ─ Т. 4, № 1. ─ P. 35-37.
21. Sinkov M.V. Research and development of algorithms for image inverse function of hypercomplex variable / M.V. Sinkov, Y.O.Kalinovsky, Yu.E.Boyarinova // Data Recording,Storage & Processing. ─ 2005. ─ Т. 7, № 1. ─ P. 32–42.